\documentclass{amsart}

\usepackage[mathscr]{eucal} 

\usepackage{pst-node}

\usepackage{tikz-cd} 
\usepackage{amsfonts}
\usepackage{amssymb}
\usepackage{amsmath}
\usepackage[
backend=biber,
style=alphabetic,
]{biblatex}

\usepackage{xcolor}

\usepackage[hidelinks]{hyperref} 
\hypersetup{
    colorlinks=true,
    linkcolor=red,
    citecolor=blue,
    filecolor=red,
   urlcolor=red}

\newcommand{\g}{\Gamma}

\def\bN{\mathbb{N}}

\newtheorem{theorem}{Theorem}[section]
\newtheorem{lemma}{Lemma}[section]
\newtheorem{corollary}{Corollary}[section]
\newtheorem{proposition}{Proposition}[section]

\theoremstyle{definition}
\newtheorem{definition}{Definition}[section]
\newtheorem{example}{Example}[section]

\title{ groups satisfying the minimal condition on centralizers are complete co-analytic}
\author{Hakkı Kağan Yumak}
\begin{document}

\maketitle

\begin{abstract}
In this paper, we prove that the set of groups satisfying the minimal condition on centralizers is a complete co-analytic subset of the Polish space of countable groups. This answers a question of Wesolek and Williams.
\end{abstract}
\section{Introduction}
Let $F_\infty$ be the free group generated by $X_\infty=\{x_i:i\in\bN\}$. We denote the set of all normal subgroups of $F_\infty$ by $\mathcal{N}(F_\infty)$. Identifying $\mathcal{P}(F_\infty)$ with $2^{F_\infty}$, and endowed with the subspace topology inherited from the product topology of 2$^{F_\infty}$, $\mathcal{N}(F_\infty)$ becomes a \textit{Polish space}, that is, a separable, completely metrizable topological space. Since each countable group is a homomorphic image of $F_\infty$, one can view $\mathcal{N}(F_\infty)$ as \textit{the space of (marked) countable groups}. For details, see \cite[Section 2.1]{MR3668055}.

Having constructed a Polish space of countable groups, one may investigate the definability of group properties within the space from a descriptive set-theoretic point of view, e.g., see \cite{MR3828881,MR3668055} and, for the space of finitely generated groups, \cite{MR4518646}.
 In this context, a group property $P$ can be considered as the subset
\[\{N \in \mathcal{N}(F_\infty):\ F_\infty/N \text{ has the property }P\}\]
of the Polish space $\mathcal{N}(F_\infty)$. In particular, group properties that define Borel sets may be thought of simpler and more explicit properties than those properties that define \textit{analytic sets}, that is, projections of Borel sets.

For example, while the set of amenable groups is Borel, the set of elementary amenable groups is non-Borel \textit{co-analytic}, that is, its complement is analytic \cite[Theorem 1.3]{MR3668055}, which also provides an alternative proof for the existence of an amenable group that is not elementary amenable. Intuitively speaking, this means that being elementary amenable cannot be defined without a universal quantification over an uncountable Polish space.

Another example result is that the set of complete groups is complete co-analytic and consequently is non-Borel \cite{MR3828881}. Here an analytic (respectively, co-analytic) set is said to be \textit{complete} if it is the most complex among analytic (respectively, co-analytic) sets in the sense of Borel reducibility.

In this paper, we continue the investigation initiated in \cite{MR3668055} and analyze the set of groups with minimal condition on centralizers.
\begin{definition}
    A group $G$ is said to satisfy the \textit{minimal condition on centralizers} if every descending chain of centralizers $(\operatorname{C}(A_i))_{i\in\bN}$ of $G$ stabilizes, where $A_i\subseteq G$.
\end{definition}
We shall denote the set of countable groups satisfying the minimal condition on centralizers by $\mathscr{M}_C$, that is,
$$\mathscr{M}_C=\{N\unlhd F_\infty: F_\infty/N \:\text{satisfies the minimal condition on centralizers}\}.
$$

It was shown in \cite[Theorem 3.12]{MR3668055} that $\mathscr{M}_C$ is co-analytic and non-Borel. Furthermore, it was asked whether $\mathscr{M}_C$ is complete \cite[Question 8.5.]{MR3668055}. The main result of this paper is the following.

\begin{theorem}\label{main theorem}
    $\mathscr{M}_C$ is a complete co-analytic subset of $\mathcal{N}(F_\infty).$
\end{theorem}

\section{Preliminaries}
\subsection{Right-Angled Artin Groups}
In this subsection, we will define right-angled Artin groups and give the necessary facts to prove the main theorem. For a more detailed investigation of basic properties of these groups, we refer the reader to \cite{MR634562}. 
\begin{definition}
    A \textit{graph} $\Gamma=(V(\Gamma), E(\Gamma))$ is defined as an ordered pair where $V(\Gamma)$ denotes a non-empty set of vertices and $E(\Gamma)\subseteq V(\Gamma)\times V(\Gamma)$ is a set of edges that forms a symmetric and irreflexive relation. 
    Given an element $u\in V(\g)$, we write $Link_\g(u)$ (or simply $Link(u)$) to denote the set of all vertices adjacent to $u$, that is, $$Link_\g(u)=\{v\in V(\g):(u,v)\in E(\g)\}.$$ 
Given $S\subseteq V(\g)$, the \textit{subgraph of $\Gamma$ induced by $S$} is defined as the graph $\Gamma_S$ with vertex set $S$ and edge set $E\cap(S\times S)$.

\end{definition}

\begin{definition}
    Let $\Gamma=(V(\Gamma),E(\Gamma))$ be a graph. Then the \textit{right-angled Artin group} of the graph $\Gamma$, is the group $G(\Gamma)$ with the following presentation:
    $$G(\Gamma)=\langle V(\Gamma)\ |\ [v,u]=e \:\text{for}\: (u,v)\in E(\Gamma)\rangle.$$
\end{definition}
Let us provide some basic examples of such groups. For example, in the case of a complete graph, the corresponding right-angled Artin group is isomorphic to the direct sum of $\mathbb{Z}$ with the index set as the set of vertices, whereas, in the case of a trivial graph, the associated right-angled Artin group is isomorphic to the free group whose rank is the cardinality of the set of vertices.

Every element in a right-angled Artin group $G(\Gamma)$ can be represented uniquely up to shuffling of adjacent vertices \cite[Lemma 2.2. (i)]{MR634562}. This representation is called a \textit{minimal form}. Given $u \in G(\Gamma)$, the \textit{length} of $u$, denoted by $|u|$, is the length of the minimal form of $u$. The \textit{support} of $u$, denoted by $\operatorname{supp}(u)$, is the set of generators occurring in the minimal form of $u$, that is,
\[
\operatorname{supp}(u)=\{v\in V(\Gamma): v\:\text{occurs in the minimal form of $u$}\}\subseteq G(\g).
\] 
The support of $u$ is well-defined by \cite[Lemma 2.2. (ii)]{MR634562}. Also, $$\operatorname{supp}(e)=\emptyset.$$ Given $A\subseteq G(\g)$ and $u\in G(\g)$, we define the \textit{support} of $A$ as
\[\operatorname{supp}(A)=\bigcup_{u \in A} \operatorname{supp}(u)\]
and the \textit{link} of $u$ as
$$\operatorname{link}(u)=\displaystyle\bigcap_{v\in \operatorname{supp}(u) }Link(v)$$
and $$\operatorname{link}(e)=V(\g)$$
It is easily seen that $Link_{\Gamma}(u)=link(u)$ for every $u \in V(\Gamma) \subseteq G(\Gamma)$.

\begin{definition} Let $G(\Gamma)$ be a right-angled Artin group. An element $u \in G(\Gamma)$ is said to be \textit{cyclicly reduced} if it is of minimal length among elements in the conjugacy class of $u$.
\end{definition}
\begin{proposition} \label{cyclyic red}\textnormal{\cite{MR1023285}} 
Let $G(\g)$ be the right-angled Artin group of a graph $\g$ and $u\in G(\Gamma)$. Then, there exists unique cyclicly reduced element $\operatorname{CR}(u)\in G(\Gamma)$ such that $u=\operatorname{CR}(u)^w  $ and $|u|=|\operatorname{CR}(u)|+2|w|$ for some $w\in G(\Gamma)$ where $$\operatorname{CR}(u)^w= w\operatorname{CR}(u)w^{-1}.$$

\end{proposition}
This unique element $\operatorname{CR}(u)$ is referred to as the \textit{cyclicly reduced element} for $u$.
\begin{definition}
  Let $G(\g)$ be a right-angled Artin group and $H \leqslant G(\g)$. The subgroup $H$ is said to be \textit{graphical} if there exists $S \subseteq V(\Gamma)$ such that $H=\langle S \rangle$. In this case, we denote $H$ by $G_S(\Gamma)$ and call it the \textit{graphical} subgroup generated by $S$.
 
    The subgroup $H$ is \textit{parabolic} if it is conjugate to a graphical subgroup.
\end{definition}
It is easily checked that a graphical subgroup $G_S(\Gamma)$ of $G(\Gamma)$ is canonically isomorphic to the right-angled Artin group of the induced subgraph $\Gamma_S$ of $\Gamma$ with the induced homomorphism of the canonical map 
\begin{align*}
    \varphi: S&\rightarrow G_S(\g)\\
    v&\mapsto v
\end{align*}
Hence, two generators of a right-angled Artin group commute if and only if they are adjacent in the graph. Moreover, given $S\subseteq V(\g)$, $S$ is finite if and only if $\langle S\rangle$ is finitely generated. To see the non-trivial direction of this equivalence, assume to the contrary that $S$ is infinite but $\langle S\rangle $ is finitely generated. Then, by the above isomorphism, $G(\g_S)$ is finitely generated. As $S$ is infinite, there is a canonical epimorphism from $G(\g_S)$ to the non-finitely generated group
$$G(K_\bN)\cong\displaystyle\bigoplus_{i\in\bN} \mathbb{Z}$$ where $K_\bN$ is the complete graph on $\bN$, which is a contradiction.

Every parabolic subgroup is also a right-angled Artin group. Thus, one can use the group-theoretic properties of right-angled Artin groups to study such subgroups. Next, we will give some basic facts regarding parabolic groups.

\begin{lemma}\textnormal{\cite[Lemma 3.3.]{MR3365774}}\label{graphint}
Let $G(\g)$ be a right-angled Artin group and $\mathcal{S}$ be a collection of subsets of $V(\g)$. Then, $$\bigcap_{S\in\mathcal{S}}G_S(\g)=G_T(\g)$$ where $T=\displaystyle\bigcap_{S\in \mathcal{S}} S.$
\end{lemma}

\begin{lemma}\label{intparabol} \textnormal{\cite[Proposition 3.4.]{MR3365774}} 
Let $G(\g)$ be a right-angled Artin group. Let $S,T\subseteq V(\g)$ and $g\in G(\g)$. Then there exists $P\subseteq S\cap T$ and $h\in G_T(\g)$ such that
\[gG_S(\g)g^{-1}\cap G_T(\g)= hG_P(\g)h^{-1}.\]
Moreover, $P=S\cap T\cap \operatorname{link}(k)$ for some $k\in G(\g)$ such that $g=hk$.
\end{lemma}

It follows from the previous lemma that the intersection of two parabolic subgroups is also a parabolic subgroup.

\begin{lemma}\label{invfix} \textnormal{\cite[Lemma 3.9.]{MR3365774}} Let $G(\g)$ be a right-angled Artin group. If $H$ is a parabolic subgroup of $G(\g)$ and $H^{g}\subseteq H$ for some $g\in G(\g)$, then $H^{g}= H$.
\end{lemma}

The main construction of this paper is based on the following theorem that characterizes centralizers of elements in a right-angled Artin group.

\begin{theorem}\label{centralizer theorem} \textnormal{\cite[The Centralizer Theorem]{MR1023285}} 
   Let $G(\g)$ be a right-angled Artin group. Let $u\in G(\g)$ be such that $u=\operatorname{CR}(u)^p$ for some $p\in G(\g)$. Then 
    $$
    \operatorname{C}(u)=\bigg(\langle A\rangle \oplus G_S(\g)\bigg)^p
    $$ where $\langle A\rangle\leq \langle\operatorname{supp(CR}(u))\rangle$  is a free abelian group of finite rank and $S=\operatorname{link(CR}(u)).$ Moreover, $\operatorname{C}(u)$ is isomorphic to a right-angled Artin group.
\end{theorem}

\begin{lemma}\label{rootsofcentralizer}\textnormal{\cite[Corollary 2.11]{MR4988587}}
    Let $G(\g)$ be a right-angled Artin group. For all $u\in G(\g)$ and $n>0$, $$
    \operatorname{C}(u)=\operatorname{C}(u^n).
    $$
    Thus, given $A\subseteq G(\g)$, $w\in G(\g)$ and $n>0$, if $w^n\in\operatorname{C}(u)$, then $w\in \operatorname{C}(A).$
\end{lemma}

\begin{theorem}\textnormal{\cite[Theorem 1.1.]{MR634562}}
    Every abelian subgroup of a right-angled Artin group is free abelian.
\end{theorem}

\begin{corollary}\label{abelianaltıstabil}
    Let $G(\g)$ be a right-angled Artin group. For all finitely generated abelian subgroups $H$ of $G(\g)$ and a subset $A$ of $G(\g)$ with $\operatorname{C}(A)\leq H$, if $\operatorname{rank}(\operatorname{C}(A))=\operatorname{rank}(H)$, then $\operatorname{C}(A)=H.$ Moreover, every descending chain $(\operatorname{C}(A_i))_{i\in\bN}$ of $G(\g)$ with $\operatorname{C}(A_i)\leq H$ for some $i\in\bN$ stabilizes.
\end{corollary}
\begin{proof}
    As $\operatorname{rank}(\operatorname{C}(A))=\operatorname{rank}(H)$, we have that $H/\operatorname{C}(A)$ is finite. Hence, for all $w\in G(\g)$, there exists $n>0$ such that $w^n\in \operatorname{C}(A)$. So, by Lemma \ref{rootsofcentralizer}, $w\in \operatorname{C}(A)$. 
\end{proof}

Before we proceed, let us provide an example that illustrates some of the notions introduced in this subsection. 

\begin{example}
   Consider the graph $\g=(V(\g),E(\g))$ where $V(\g)=\{x_1,x_2,x_3,x_4\}$ and $E(\g)=\{(x_i,x_j), (x_3,x_4),(x_4,x_3)\in V(\g)^2:i\neq j\leq 3\}.$ Then, the right-angled Artin group $G(\g)$ is isomorphic to $\mathbb{Z}^2\ast_{\mathbb{Z}}\mathbb{Z}^3$. A minimal form of the word $$u=x_1^{-1}x_3x_4x_2x_1x_3^{-1}$$ is  $x_1^{-1}x_4x_2x_1$. Also, $|u|=4$, $\operatorname{supp}(u)=\{x_1,x_2,x_4\}$,  and $\operatorname{link}(u)=\{x_3\}.$ The word $u$ is not cyclicly reduced since $|u^{x_1}|=|x_4x_2|=2$. Also, $\operatorname{CR}(u)=x_4x_2$ and $u=\operatorname{CR}(u)^{x_1^{-1}}$. By the Centralizer Theorem, $$\operatorname{C}(u)=\left(\langle x_4x_2\rangle\oplus \langle x_3\rangle\right)^{x_1^{-1}}.$$
\end{example}

\subsection{Some facts from descriptive set theory}
In this paper, we will provide some necessary facts from descriptive set theory. For a general reference, we refer the reader to \cite{MR1321597}.

\begin{definition}
    Let $X, Y$ be Polish spaces, $A\subseteq X$ and $B\subseteq Y$. A Borel map $f:X\rightarrow Y$ is called \textit{Borel reduction} of $A$ to $B$ if $f^{-1}(B)=A.$ In this case, we say that $B$ reduces $A$.
\end{definition}
\begin{definition}
    Let $B$ be a subset of a Polish space $Y$. The set $B$ is called \textit{complete analytic} (respectively, \textit{co-analytic}) if $B$ is an analytic (respectively, co-analytic) subset and for every standard Borel space $X$ and every analytic (respectively, co-analytic) subset $A$ of $X$, there is a Borel reduction of $A$ from $B$.
\end{definition}
Moreover,  the complement of a complete analytic set is itself complete co-analytic.
Note that an analytic set that reduces a complete analytic set is itself complete analytic. In our main construction, we will use the following well-known complete analytic subset, namely, the set of ill-founded trees.

\begin{definition} Let $\bN^{<\bN}$ denote the set of finite sequences over $\mathbb{N}$. A subset $T$ of $\bN^{<\bN}$ is called a \textit{descriptive set-theoretic tree } if it is closed under taking initial segments. The tree $T$ is called \textit{ill-founded} if there exists $\alpha\in\bN^\bN$ such that $$\alpha|_n=(\alpha_0,\dots,\alpha_{n-1})\in T$$ for all $n\in\bN$. Such $\alpha$ is said to be an \textit{infinite branch} of $T$.
\end{definition}
Endowing the set of descriptive set-theoretic trees with the subspace topology induced from the product topology of $2^{\bN^{<\bN}}$ where $\mathcal{P}\left(\bN^{<\bN}\right)$ and $2^{\bN^{<\bN}}$ are naturally identified, we obtain the Polish space $Tr$ of descriptive set-theoretic trees.
\begin{theorem}\textnormal{\cite[(27.1) Theorem.]{MR1321597}}
    The set \textnormal{\textbf{IF}} of ill-founded descriptive set-theoretic trees is complete analytic in $Tr$.
\end{theorem}
\section{Reduction}\label{s3}
In order to prove that $\mathscr{M}_C$ is complete co-analytic, we will show that $\mathcal{N}(F_\infty)\setminus\mathscr{M}_C$ is complete analytic by constructing a Borel reduction of \textbf{IF} to $\mathcal{N}(F_\infty)\setminus\mathscr{M}_C$.

Let us first enumerate $X_{\infty}$ by the set $\bN^{<\bN}$ as $X_{\infty}=\{x_s\}_{s \in \bN^{<\bN}}$. Let $T$ be a descriptive set-theoretic tree. Define the normal subgroup $N(T)$ of $F_{\infty}$ associated with $T$ as follows:

\[N(T)=\big\langle\big\langle\{x_s:s\not\in T\} \cup \{[x_s,x_t]: (s,t) \in J\}\big\rangle\big\rangle\]
where $J$ is the set of $(s,t)$ such that either
\begin{itemize}
    \item $s \neq t$ are compatible and both of $s$ and $t$ are of even length, or
    \item $s$ is of even length, $t$ is of odd length and $s$ is an initial segment of $t$.
\end{itemize}
Consider the map $\varphi:Tr\rightarrow\mathcal{N}(F_\infty)$ defined by $\varphi(T)=N(T)$. We shall show that $\varphi$ is a Borel reduction of $\mathbf{IF}$ to $\mathcal{N}(F_\infty)\setminus\mathscr{M}_C$, that is, $\varphi$ is Borel and
\[
T\in\textbf{IF}\Leftrightarrow \varphi(T)\not\in\mathscr{M}_C
\]
We begin by showing that $\varphi$ is indeed continuous.

\begin{proposition}\label{phiisborel} $\varphi$ is continuous. 
\end{proposition}
\begin{proof}
 Let $(T_n)_{n\in\bN} $ be a sequence of descriptive set-theoretic trees. Suppose that the sequence converges to $T$. 
    Let $w(x_{s_1},\ldots,x_{s_k})\in F_\infty.$
Since $(T_n)_{n\in\bN}$ converges to $T$, there exists $M\in\bN$ such that for all $n\geq N$ and all $i\leq k$,  $$
    s_i\in T_n \Leftrightarrow s_i\in T.
    $$ 
    So, for all $n\in\bN$, the subtree
    $$S_n:=\{s: s\subseteq s_i\in T_n \: \text{for some } i\leq k\}$$ is equal to the subtree
    $$S:=\{s: s\subseteq s_i\in T \: \text{for some } i\leq k\}.$$
    Hence, $N(S_n)=N(S)$ for all $n\geq M.$

  Suppose that $w(x_{s_1},\ldots,x_{s_k})\in N(T).$ Assume on the contrary that there exists $n\geq M$ such that $w(x_{s_1},\ldots,x_{s_k})\not\in N(T_n).$ Consider the following homomorphism 
 \begin{align*}
       \psi:F_\infty&\rightarrow F_\infty/N(T_n)\\
      x_s&\mapsto \begin{cases}
           x_sN(T_n),\quad &\text{if}\: s\in S_n\\
           e,\quad&\text{if} \:s\not\in S_n
       \end{cases}
   \end{align*}
Then, it is easy to observe that  $\ker(\psi)=N(S_n)$.
As $N(T)\subseteq N(S)=N(S_n)$, we have that $\psi(w(x_{s_1},\ldots,x_{s_k}))=e.$ So, $w(x_{s_1},\ldots,x_{s_k})\in N(T_n)$, which is a contradiction. Thus, $w(x_{s_1},\ldots,x_{s_k})\in N(T_n)$ for all $n\geq M.$

    Suppose that $w(x_{s_1},\ldots,x_{s_k})\not\in N(T).$ 
   Assume on the contrary that there exists $n\geq M$ such that $w(x_{s_1},\ldots,x_{s_k})\in N(T_n)$. Since $N(T_n)\subseteq N(S_n)=N(S)$, $w(x_{s_1},\ldots,x_{s_k})\in N(S_n).$ 
   Now, define the homomorphism 
   \begin{align*}
       \Psi:F_\infty&\rightarrow F_\infty/N(T)\\
      x_s&\mapsto \begin{cases}
           x_sN(T),\quad &\text{if}\: s\in S\\
           e,\quad&\text{if} \:s\not\in S
       \end{cases}
   \end{align*}
  Since $\ker(\Psi)=N(S),$ we have $\Psi(w(x_{s_1},\ldots,x_{s_k}))=e$. So, $w(x_{s_1},\ldots,x_{s_k})\in N(T)$, which is a contradiction. Therefore, $w(x_{s_1},\ldots,x_{s_k})\not\in T_n$ for all $n\geq M$. Hence, $\varphi$ is continuous.       
\end{proof}

In order to show that $\varphi$ is a reduction, we shall need a series of lemmas and corollaries. First, observe that $G(T)=F_\infty/N(T)$ is the right-angled Artin group of the graph $$\g(T)=(V(\g(T)),E(\g(T)))$$ where $V(\g(T))=\{x_s: s \in T\}$ and $E(\g(T))=\{(x_s,x_t),(x_t,x_s): (s,t)\in J\}$.
This observation allows us to use properties of right-angled Artin groups.

\begin{lemma}\label{oddincomp}
    Let $u\in G(T)$.  If $\operatorname{supp(CR}(u))$ contains 
 $x_s,x_t$ for some incompatible $s,t\in T$, or
 $x_t$ for some $t\in T$ of odd length,
    then $\operatorname{C}(u)$ is finitely generated abelian subgroup.
\end{lemma}
\begin{proof}
    As $\operatorname{C}(u)$ and $\operatorname{C(CR(}u))$ are conjugate, it is enough to show that $\operatorname{C(CR(}u))$ satisfies the condition.
    
    According to Theorem \ref{centralizer theorem}, the centralizer $\operatorname{C(CR(}u))$ of $\operatorname{CR}(u)$ is of the form
    $$\langle A\rangle \oplus \langle\operatorname{link(CR(}u))\rangle$$
  where $\langle A\rangle$ is a finite rank free abelian subgroup. To deduce the result, we will show that $\langle\operatorname{link(CR(}u))\rangle$ is finitely generated abelian subgroup.
 
\underline{Case 1:} Suppose that there exists $x_t\in \operatorname{supp(CR}(u))$ with $t$ of odd length. Then,  by the construction, we have that  $$
\operatorname{link}(x_t)=\{x_s\in X_\infty: s\in T,\  s\subseteq t\:\text{and}\: s\:\text{ is of even length}\}
    $$
    is finite. We also have $\operatorname{link(CR}(u))$ is finite as  $\operatorname{link(CR}(u)) \subseteq\operatorname{link}(x_t).$ 
Therefore, $\operatorname{C}(u)$ is finitely generated.
Now, we will show that $
\langle A\rangle\oplus \langle \operatorname{link}(x_t)\rangle
$ is abelian since $\operatorname{C}(\operatorname{CR}(u))$ contained in this subgroup. To deduce that, it suffices to show that $\langle \operatorname{link}(x_t)\rangle$ is abelian.

Let $x_s, x_k \in \operatorname{link}(x_t)$. By construction, $s,k\subseteq t$ and both have even length. So, $(s, k) \in J$, and thus $x_s$ and $x_k$ commute. Consequently, $\langle \operatorname{link}(x_t) \rangle$ is abelian.

\underline{Case 2:}
Suppose that there exists $x_s,x_t\in\operatorname{supp(CR}(u))$ such that $t$ and $s$ are incompatible.
Observe that
$$
\operatorname{link}(\operatorname{CR}(u)) \subseteq \operatorname{link}(\{x_s, x_t\}) = \{x_k : k \subseteq s, t \text{ and } |k| \text{ is even}\}.
$$
Here, $\langle \{x_k : k \subseteq s, t \text{ and } |k| \text{ is even}\} \rangle$ is a finitely generated abelian subgroup by the construction. Therefore, $\operatorname{C}(\operatorname{CR}(u))$ is also finitely generated and abelian, and thus $\operatorname{C}(u)$ is a finitely generated abelian subgroup as well.
\end{proof}
Throughout the rest of paper, for each $t \in T$, let $X_t$ denote the set
$$
\{x_k : t \subseteq k,\; k \in T\} \cup \{x_k : k \subseteq t\; \text{and}\; |k|\; \text{is even}\}\subseteq V(\g(T)).
$$
Additionally, for each $u \in G(T)$, let $t_u$ denote the maximum element (if it exists) in
$
\{s : x_s \in \operatorname{supp}(\operatorname{CR}(u))\}.
$
\begin{lemma}\label{charachcentrza}
    Let $u\in G(T)$ with $u=\operatorname{CR}(u)^{p}$ for some $p\in G(T)$. If there exists $r\in T$ with
$\operatorname{supp(CR}(u))\subseteq\{x_{s}\in G(T):s\subseteq r \:\text{and}\: |s|\:\text{is even} \},$ then $\operatorname{C}(u)=\langle X_{t_u}\rangle^{p}$.
\end{lemma}
\begin{proof}
According to  Theorem \ref{centralizer theorem} and the construction, the centralizer $\operatorname{C(CR}(u))$ is of the form
    $$
\left(\bigoplus_{x_s\in \operatorname{supp(CR}(u))}\langle x_s\rangle \right) \oplus \langle \operatorname{link(CR}(u))\rangle 
    $$
As $\operatorname{supp(CR}(u))\subseteq\{x_{s}\in G(T):s\subseteq r \:\text{and}\: |s|\:\text{is even}\}$ and by the construction, we have that $\operatorname{link(CR}(u))=X_{t_u}\setminus\operatorname{supp(CR}(u))$. 
 Therefore, $\operatorname{C(CR}(u))=\langle X_{t_u}\rangle$.
\end{proof}

\begin{lemma}\label{incompfgabel}
    Let $u,w\in G(T)$ with $u=\operatorname{CR}(u)^{p}$ and $w=\operatorname{CR}(w)^{q}$ for some $p,q$. If there exists $x_s\in\operatorname{supp(CR}(u))$ and  $x_t\in\operatorname{supp(CR}(w))$ with $s$ and $ t$ are incompatible, then  $\operatorname{C}(\{u,w\})$ is finitely generated abelian group.
\end{lemma}
\begin{proof}
     Suppose that there exists incompatible $s,t\in T$ such that  $x_s\in \operatorname{supp(CR}(u))$ and $x_t\in \operatorname{supp(CR}(w))$. As $\operatorname{C}(\{u,w\})\leq \operatorname{C}(u),\operatorname{C}(w)$ and by Lemma \ref{oddincomp}, suppose without loss of generality, there exists $r,k\in T$ such that 
     \begin{align*}
         \operatorname{supp(CR}(u))\subseteq \{x_{s}&\in G(T):|s|\:\text{is even and}\: s\subseteq r\}\\
         &\text{and}\\
         \operatorname{supp(CR}(w))\subseteq \{x_{s}&\in G(T):|s|\:\text{is even and}\: s\subseteq k\}.
     \end{align*} Thus, by Lemma \ref{charachcentrza}, we have that $$
    \operatorname{C}(u)=\langle X_{t_u}\rangle^{p}\:\text{and}\:\operatorname{C}(w)=\langle X_{t_w}\rangle^{q}.
    $$  
    So,
    by Lemma \ref{intparabol}, there exists $P\subseteq X_{t_u}\cap X_{t_w} $ and $h\in G(T)$ such that  $$\operatorname{C}(\{u,w\})=G_P(\g(T))^h.$$ 
As $G_P(\g(T))^h\leq \langle X_{t_u}\cap X_{t_w} \rangle^h$, it is enough to check whether $\langle X_{t_u}\cap X_{t_w} \rangle$ is a finitely generated abelian subgroup to deduce the result. Since there exist incompatible elements $s\in\{l : x_l \in \operatorname{supp}(\operatorname{CR}(u))\}$ and $t\in\{l : x_l \in \operatorname{supp}(\operatorname{CR}(w))\}$,  and $t_u$ and $t_w$ are maximum elements in these sets, respectively, it follows that $t_u$ and $t_w$ are also incompatible.
So, $$X_{t_u}\cap X_{t_w}=\{x_k: k\subseteq t_u,t_w\:\text{and}\: |k|\:\text{is even} \}.$$ Then, by the construction, $\langle X_{t_u}\cap X_{t_w} \rangle$ is finitely generated abelian subgroup.
\end{proof}

It follows that if a descending chain of centralizers does not stabilize, the supports of its elements must lie along a single branch and consist of even-length nodes. That is, by combining Lemma \ref{incompfgabel}, Lemma \ref{oddincomp} and Corollary \ref{abelianaltıstabil}, we have the following: 
\begin{corollary}\label{compandeven}
Let $(\operatorname{C}(A_i))_{i\in\bN} $ be a descending chain of centralizers in $G(T)$.
If the set $$\displaystyle\bigcup_{i\in\bN}\bigcup_{w\in A_i}\operatorname{supp(CR}(w))$$ contains
\begin{itemize}
    \item $x_s,x_t$ for some incompatible $s,t\in T$, or
    \item  $x_t$ for some $t\in T$ of odd length,
\end{itemize}
then the chain stabilizes.
\end{corollary}

\begin{lemma}\label{stabfinite}
Let $(\operatorname{C}(A_i))_{i\in\bN} $ be a descending chain of centralizers in $G(T)$. If the set $$\displaystyle\bigcup_{i\in\bN}\bigcup_{w\in A_i}\operatorname{supp(CR}(w))$$ is finite, then the chain stabilizes.
\end{lemma}
\begin{proof}
Assume on the contrary that the chain $(\operatorname{C}(A_i))_{i\in\bN} $ does not stabilize. Since the union is finite and the chain does not stabilize,  by the previous corollary,  we have that $$\displaystyle\bigcup_{i\in\bN}\bigcup_{w\in A_i}\operatorname{supp(CR}(w))\subseteq X_0= \{x_{s}\in G(T): s\subseteq r\: \text{and} \:|s|\:\text{is even}\: \}$$  for some $r\in T$.

Let $i\in\bN$ and  $u\in A_i$ and write $u=\operatorname{CR}(u)^{p_u}$ for some $p_u\in G(T)$. 
Since $X_0$ forms a finite chain and $\operatorname{supp(CR}(u))\subseteq X_0$,  the maximum element $t_u\in T$ of 
$
\{s : x_s \in \operatorname{supp}(\operatorname{CR}(u))\}
$ exists. Also,  we have that  $$\operatorname{C}(u)=\langle X_{t_u}\rangle^{p_u}\: \text{and} \:
\langle X_{t_u}\rangle =\operatorname{C}(x_{t_u})$$ by Lemma \ref{charachcentrza}.
%Hence, for each $u\in A_i$, the centralizer of $u$ is conjugate to the centralizer of $x_{t_u}$, where ${t_u}$ is the maximum element of $\{t\in T: x_t\in \operatorname{supp(CR}(u))\}.$ In other words, $$\operatorname{C}(u)=\operatorname{C}(x_{t_u})^{p_u}=\langle X_{t_u}\rangle^{p_u}$$ where $X_t\subseteq V(\g(T)).$
    Therefore, the centralizer of $A_i$ can be expressed as an intersection of conjugates of some elements in $X_0.$ That is, $$\operatorname{C}(A_i)=\displaystyle\bigcap_{w\in A_i}\operatorname{C}(x_{t_w})^{p_w}$$ for some $x_{t_w}\in X_0$ and $p_w\in G(T).$

We aim to show that the centralizer of $A_i$ equals the conjugate of the centralizer of $x_s\in X_0$. To do that, we first analyze the intersection of two distinct conjugates of the centralizer of the same element. 
We claim that for $x_t \in X_0$ and $p, q \in G(T)$, if $\operatorname{C}(A_i)$ is contained in the intersection $\operatorname{C}(x_t)^p \cap \operatorname{C}(x_t)^q$, then $\operatorname{C}(x_t)^p = \operatorname{C}(x_t)^q$.

Suppose that $\operatorname{C}(A_i)$ is contained in the intersection $\operatorname{C}(x_t)^p \cap \operatorname{C}(x_t)^q$. Also, we have that $$
\operatorname{C}(x_t)^{p}=\langle X_{t}\rangle^{p} \:\text{and}\: \operatorname{C}(x_t)^{q}=\langle X_{t}\rangle^{q}
$$
So, $$\left(\operatorname{C}(x_t)^{q^{-1}p}\cap \operatorname{C}(x_t)\right)^{q}=\operatorname{C}(x_t)^{p}\cap \operatorname{C}(x_t)^{q}$$
Then, by Lemma \ref{intparabol}, there exists $P\subseteq X_t$ and $h\in \langle X_t\rangle$ such that $$\operatorname{C}(x_t)^{q^{-1}p}\cap \operatorname{C}(x_t)=G_{P}(\g(T))^h$$ 
Moreover, we have  $P=X_t\cap \operatorname{link}(g)$ where $hg=q^{-1}p$ for some $g\in G(T)$. 
  It follows from the construction that if $\operatorname{supp}(g)$ contains $x_k$ for some $k\in T$ of odd length or incompatible with $t$,  then $P\subseteq \{x_s: s\subseteq l\:\text{and}\: |s|\:\text{is even}\}$ for some $l\in T$ in any case. Then, by construction, $G_{P}(\g(T))^h$ is a finitely generated abelian subgroup. So, by Corollary \ref{abelianaltıstabil}, the chain $(\operatorname{C}(A_i))_{i\in\bN}$ stabilizes, which is a contradiction.
  Therefore,  $$\operatorname{supp}(g)
   \subseteq X_t.$$ 
   So, $hg=q^{-1}p\in \langle X_t\rangle=\operatorname{C}(x_t)$, which implies   $\operatorname{C}(x_t)^{p}=\operatorname{C}(x_t)^{q}$. This proves the claim.
 
Thus, the centralizer of $A_i$ is of the form $$
\bigcap_{x_t\in T_i} \operatorname{C}(x_t)^{p_t}
$$ 
 for some $p_t\in G(T)$ where $
T_i=\{x_t: t=t_u\: \text{for some } u\in A_i\}\subseteq X_0.$

Next, we will investigate the conjugates of centralizers of different elements. 

Let $p_t,p_s\in G(T)$ with $x_t\neq x_s\in T_i$. As $x_t,x_s\in X_0$, either $s\subseteq t$ or $t\subseteq s$. Without loss of generality, suppose  $s\subseteq t$. Then $\operatorname{C}(x_t)\subseteq \operatorname{C}(x_s)$. So,
$$
\operatorname{C}(x_t)^{p_t}\cap \operatorname{C}(x_s)^{p_s}\subseteq \operatorname{C}(x_s)^{p_t}\cap \operatorname{C}(x_s)^{p_s}.
$$ Our previous claim shows that $\ \operatorname{C}(x_s)^{p_t}=\operatorname{C}(x_s)^{p_s}.$ 
Therefore, we conclude that $$
\operatorname{C}(x_t)^{p_t}\cap \operatorname{C}(x_s)^{p_s}=(\operatorname{C}(x_t)\cap \operatorname{C}(x_s))^{p_t}=\operatorname{C}(x_t)^{p_t}.
$$

Hence, we obtain the following:
$$
\operatorname{C}(A_i)=\operatorname{C}(x_{t_i})^{p_i}
$$ for some  $p_i\in G(T)$ and $x_{t_i}\in X_0$  where $t_i$ is the maximum element in $\{t:x_t\in T_i\}.$

As $X_0$ is finite, there exists $x_t \in X_0$ such that $\operatorname{C}(A_i) = \operatorname{C}(x_t)^{p_i}$ for infinitely many $i \in \mathbb{N}$.
As $$\operatorname{C}(x_{t})^{p_{i_n}}\supseteq \left(\operatorname{C}(x_{t})^{p_{i_n}}\right)^{{p_{i_{n+1}}p_{i_n}^{-1}}}=\operatorname{C}(x_{t})^{p_{i_{n+1}}}$$ for all $n\in\bN$ and by Lemma \ref{invfix}, we have that $$\operatorname{C}(x_{t})^{p_{i_n}}= \operatorname{C}(x_{t})^{p_{i_{n+1}}}.$$ Since $(\operatorname{C}(x_{t})^{p_{i_n}})_{n\in\bN} $ is a stabilizing subchain of the descending chain  $(\operatorname{C}(A_i))_{i\in\bN}$, it follows that $(\operatorname{C}(A_i))_{i\in\bN}$ also stabilizes, which is a contradiction.
\end{proof}

\begin{theorem}\label{phiisreduction} $\varphi$ is a reduction from $\mathbf{IF}$ to $\mathcal{N}(F_\infty)\setminus\mathscr{M}_C$. 
\end{theorem}
\begin{proof}
    Suppose  $T$ is ill-founded and $\alpha\in[T]$. For each $i\in\bN$, define $$A_i=\{x_{\alpha|_{2k}}\in G(T): k\leq i \}.$$
    By Lemma \ref{charachcentrza}, we have that 
    $\operatorname{C}(x_{\alpha|_{2k}})=\langle X_{\alpha|_{2k}}\rangle$ for all $k\in\bN$ where  
    $$
    X_{\alpha|_{2k}}=\{x_s: \alpha_{|2k}\subseteq s\}\cup\{x_{\alpha|_{2n}}: n\leq k\}.
    $$
    Therefore, by Lemma \ref{graphint},$$\operatorname{C}(A_i)=\langle\{x_s: s\supseteq \alpha_{|2i}\}\cup\{x_{\alpha|_{2n}}: n\leq i\}\rangle.$$ Since $x_{\alpha|_{2i+1}}\in \operatorname{C}(A_i)\setminus \operatorname{C}(A_{i+1})$, we obtain that $(\operatorname{C}(A_i))_{i\in\bN}$ is a strictly descending chain of centralizers. So, $\varphi(T)\not\in\mathscr{M}_C.$

    Conversely, suppose there exists a descending chain $(\operatorname{C}(A_i))_{i\in\bN} $ of centralizers that does not stabilize. Any two elements $s,t\in T$ with $x_s,x_t$ are contained in $$\displaystyle\bigcup_{i\in\bN}\bigcup_{w\in A_i}\operatorname{supp(CR}(w))$$ are  compatible by Corollary \ref{compandeven}. Moreover, this union 
    is infinite by Lemma \ref{stabfinite}. Therefore, one can recursively construct an infinite branch of $T$.
\end{proof}

Since $\mathbf{IF}$ is complete analytic, it now follows from Proposition \ref{phiisborel} and Theorem \ref{phiisreduction} that $\mathcal{N}(F_\infty)\setminus\mathscr{M}_C$ is complete analytic as well. Thus, $\mathscr{M}_C$ is complete co-analytic.

\subsection*{Acknowledgments} I would like to thank Gökhan Benli and Kaan Doğanay for helpful mathematical discussions. I also thank Burak Kaya for his feedback on an earlier version of this manuscript.
\printbibliography
\end{document}